\documentclass{article}
\usepackage[utf8]{inputenc}
\usepackage{amsmath,amsthm,amssymb,amsfonts}
\usepackage[usenames]{color}

\usepackage{graphicx}
\usepackage{todonotes}
\usepackage{tikz} 
\usepackage{comment}

\newtheorem{thm}{Theorem}

\def\1{{\bf 1}}

\title{Poisson approximation by coupling
}
\begin{document}

\maketitle

Rinaldo B. Schinazi
\footnote{Department of Mathematics, University of Colorado, Colorado Springs, CO 80933-7150, USA.
E-mail: rinaldo.schinazi@uccs.edu}

\begin{abstract}
    It is well known that a binomial $(n,p)$ can be approximated by a Poisson distribution with parameter $np$. The typical approach in undergraduate probability texts is to show a convergence result for the distribution of the binomial as $n$ goes to infinity and $np$ converges to some $\lambda$. In this note we use instead the coupling technique to show a much more general result. Moreover, we only use elementary results from probability.
    
\end{abstract}

\section{Introduction}

It is well known that a binomial random variable $B$ with parameters $(n,p)$ can be approximated by a Poisson random variable $L$ with parameter $\lambda=np$.  To show this result, undergraduate texts typically compute the limit of $P(B=k)$ as $n\to\infty$ and $np\to\lambda$ and show that this limit is $P(L=k)$ for any fixed integer $k$, see for instance \cite{Durrett} or \cite{Pitman}.  This shows the convergence in distribution of a binomial to a Poisson distribution. But this is not really an approximation result as it does not provide a bound on the error. By using the coupling technique one gets a much more general result and in particular a bound on the error valid for all $n$ and $p$. More precisely, we prove the following theorem.

\begin{thm}
    Let $B$ be a binomial random variable with parameters $(n,p)$ and $L$ be a Poisson random variable with parameter $np$. Let $D$ be a set of positive integers. Then, for all $n\geq 1$ and $p\in (0,1)$,
$$\left|P(L\in D)-P(B\in D)\right|\leq np^2. $$
\end{thm}

Since $D$ is an arbitrary set of integers and the bound $np^2$ is uniform in $D$ the maximum error by approximating the distribution of a binomial $(n,p)$ by a Poisson $np$ is smaller than $np^2$. 

Coupling is a very powerful and beautiful technique discovered by Wolfang Doeblin in the 1930's, see Appendix 4 in \cite{Lindvall}. However, for some reason coupling is mostly covered in graduate texts and research monographs, see for instance \cite{Roch}. 

The coupling we discuss in this note is not new, see for instance page 5 in \cite{Lindvall}.There, however, it is a minor point in an advanced monograph. In this note we will streamline the proof and provide all the details at an elementary level. 

It is easy to adapt the argument below to the more general case of approaching the sum $B=\sum_{i=1}^n B_i$ of independent $p_i$ Bernoulli random variables by a Poisson random variable $L$ with parameter $\lambda=\sum_{i=1}^n p_i$, see for instance \cite{Pollett}. The result is then
$$\left|P(L\in D)-P(B\in D)\right|\leq \sum_{i=1}^n p_i^2.$$
The bound above can be improved using more involved techniques. The improved bound is typically of the type $f(\lambda)\sum_{i=1}^n p_i^2$, for a certain function $f$. See for instance \cite{Barbour}.

\section{Coupling a Poisson and a Bernoulli}

Let $U$ be uniform on $(0,1)$. Let $p\in (0,1)$. For $k\geq 0$, let
$$r_k=e^{-p}\frac{p^k}{k!}\mbox { and }s_k=r_0+r_1+\dots+r_k.$$
Define the random variables $B$ and $L$ by,
$$
B=\begin{cases}
    0\mbox{ if }U<1-p\\
    1\mbox{ if }U>1-p
\end{cases}
\mbox{ and }
L=\begin{cases}
    0\mbox{ if }U<r_0\\
    k\geq 1\mbox{ if } s_{k-1}<U<s_k.
\end{cases}
$$

Hence, $B$ is a Bernoulli with parameter $p$ and $L$ is a Poisson random variable with parameter $p$. 
 This construction yields the following coupling,
$$
(L,B)=\begin{cases}
    (0,0)\mbox{ if }U<1-p\\
    (0,1) \mbox{ if }1-p<U<r_0\\
    (k,1) \mbox{ for } k\geq 1 \mbox{ if }s_{k-1}<U<s_k.
\end{cases}
$$
Therefore,
\begin{align*}
E|L-B|=&r_0-(1-p)+\sum_{k\geq 1}(k-1)r_k\\
=&r_0-(1-p)+E(L)-P(L\geq 1)\\
=&2(r_0-(1-p))\\
=&2(e^{-p}-(1-p)).\\
\end{align*}
Note that 
$$2(e^{-p}-(1-p))\leq p^2.$$
Therefore, 
$$E|L-B|\leq p^2.\leqno (1)$$

\section{Coupling a Poisson and a Binomial}

Consider a sequence $(U_i)_{1\leq i\leq n}$ of independent uniform random variables on $(0,1)$. For every $i$ in $\{1,\dots,n\}$ we use the construction in the previous section with $U_i$ (instead of $U$) to get  the pair $(L_i,B_i)$ where $L_i$ is a Poisson and $B_i$ is a Bernoulli, both with parameter $p$. Note that the random variables $(B_i)$ are mutually independent and so are the random variables $(L_i)$. 

Let $B=\sum_{i=1}^n B_i$, then $B$ is a binomial random variable with parameters $(n,p)$.
Let $L=\sum_{i=1}^n L_i$ then $L$ is a Poisson random variable with parameter $np$.

\begin{align*}
E\left|L-B\right|=&E\left |\sum_{i=1}^n (L_i-B_i)\right |\\
\leq & E\left(\sum_{i=1}^n |L_i-B_i|\right)\\
    \leq & np^2,
\end{align*}
where the last inequality comes from (1).

\section{The Poisson approximation to a binomial}

An elementary computation shows that for any random variables $B$ and $L$,
$$|P(L\in D)-P(B\in D)|\leq P(L\not = B),$$
where $D$ is a set of positive integers, see for instance Section 3, Chapter 21 in \cite{Schinazi}. Thus,
\begin{align*}
    |P(L\in D)-P(B\in D)|&\leq P(L\not= B)\\
    &=P(|L-B|\geq 1)\\
    &\leq E|L-B|\\
    &\leq np^2.
\end{align*}

\bibliographystyle{amsplain}

\end{document}